\documentclass[12pt,a4paper]{article}
\usepackage[cp1251]{inputenc}
\usepackage[T2A]{fontenc}
\usepackage{amssymb}
\usepackage{amsmath}
\usepackage{amsfonts}
\usepackage[russian]{babel}
\usepackage{amssymb, amsmath, latexsym}
\usepackage{geometry}
\begin{document}

\begin{center}

\textbf{ Sharp quadrature formulas and Nikol'skii type
inequalities\\ for rational functions}

\textbf{V.I. Danchenko and L.A. Semin}\\
vdanch2012@yandex.ru, semin.lev@gmail.com
\end{center}

\begin{quote}

\textbf{Abstract.} Sharp quadrature formulas for integrals of complex rational
functions on circles, real axis and its segments are obtained. We
also find sharp quadrature formulas for calculation of $L_2$-norms
of rational functions on such sets. Basing on quadrature formulas
for rational functions, in particular, for simple partial
fractions and polynomials, we derive sharp inequalities for
different metrics (Nikol'skii type inequalities).

\textbf{Bibliography:} 36 titles.


\end{quote}

\noindent УДК 517.5, 517.53

\begin{center}

\textbf{Точные квадратурные формулы  и неравенства \\разных метрик
для рациональных функций}

\footnote{Работа выполнена при финансовой поддержке ДРПННиТ №
1.1348.2011, РФФИ (проект № 14-01-00510), Минобрнауки России
(задание № 2014/13, код проекта: 3037). Работа выполнена в рамках
Государственного задания ВлГУ № 2014/13 в сфере научной
деятельности (тема 2868).}
\\

\textbf{В.И. Данченко, Л.А. Семин}\\
vdanch2012@yandex.ru, semin.lev@gmail.com
\end{center}

\begin{quote}

\textbf{Аннотация.} Получены точные квадратурные формулы для интегралов от
рациональных функций комплексного переменного на окружностях,
вещественной оси и ее сегментах. Найдены точные квадратурные
формулы для подсчета $L_2$-норм рациональных функций на таких
множествах. На основании квадратурных формул для рациональных
функций и, в частности, для наипростейших дробей и многочленов
получены точные неравенства разных метрик (типа неравенств С.М.
Никольского).

\textbf{Библиография:} 36 наименований.

\textbf{Ключевые слова:} точные квадратурные формулы для рациональных
функций, наипростейшая дробь, неравенства разных метрик (типа
неравенств С.М. Никольского).

\end{quote}


\bigskip
 {\bf\large 1. Формулировка основных результатов}
\medskip

{\bf 1.1.} Пусть $r>0$ и $\gamma_r:=\{z:\,|z|=r\}$. Рассмотрим
конечное множество на $\overline {\mathbb C}$ вида ${\cal Z}_n\cup
{\cal Z}_n^*$, $n\in \mathbb N$, где
$$
{\cal Z}_n:=\{z_1,\ldots,z_{n}\},\quad {\cal
Z}_n^*=\{z_1^{*},\ldots,z_{n}^{*}\}=
\{r^2/\overline{z_1},\ldots,r^2/\overline{z_{n}}\},\quad
z_1=0,\quad z_1^{*}=\infty,
   \eqno{(1)}
$$
причем множество ${\cal Z}_n$ лежит в круге $|z|<r$ (множество
${\cal Z}_n^*$ симметрично для ${\cal Z}_n$ относительно
окружности $\gamma_r$). Считаем точки $z_k$ попарно различными.
При $z\in \mathbb C$ положим
$$
B(z)=r^n\prod_{k=1}^{n}\frac{z-z_k}{r^2-z\overline{z_k}},\qquad
\mu(z)=z\frac{B'(z)}{B(z)},\quad
\mu(\zeta)=\sum_{k=1}^n\frac{r^2-|z_k|^2}{|\zeta-z_k|^2},\quad
\zeta\in\gamma_r.
  \eqno{(2)}
$$
Важную роль при построении квадратурных узлов играют корни
$\zeta_k(s,\varphi)$ уравнения $B^s(\zeta)=e^{i\varphi}$
относительно $\zeta$ при фиксированных $s\in \mathbb{N}$ и
$\varphi\in \mathbb{R}$. Будем иногда называть их {\it точками
насечки}. Очевидно, при любых $s$ и $\varphi$ имеется ровно $s\,n$
различных точек насечки и все они лежат на $\gamma_r$. Кроме того,
каждая точка  $\zeta_k(s,\varphi)$ непрерывно пробегает окружность
$\gamma_r$ при фиксированном $s$ и непрерывном изменении параметра
$\varphi$ от $0$ до $2\pi n s$. Справедлива

\medskip
{\bf Теорема 1.} {\it Пусть $s\in \mathbb{N}$ и $R$
--- рациональная функция, все полюсы которой лежат на
${\cal Z}_n\cup {\cal Z}_n^*$, причем отличные от $z_1=0$ и
$z_1^{*}=\infty$ полюсы имеют кратности не выше $s$, а в точках
$z_1$ и $z_1^{*}$ могут быть полюсы кратности не выше $s-1$. Тогда
при всех $\varphi\in \mathbb{R}$ для точек насечки
$\zeta_k=\zeta_k(s,\varphi)$, удовлетворяющих уравнению
$B^{s}(\zeta)=e^{i\varphi}$, имеем
$$
\int_{|\zeta|=r}R(\zeta)|d\zeta|=\frac{2\pi r}{s}
\sum_{k=1}^{s\,n}\frac{R(\zeta_k)}{\mu(\zeta_k)}.
 \eqno{(3)}
$$
}

\medskip
{\bf Теорема 2.} {\it В условиях теоремы $1$ для точек насечки
$\zeta_k=\zeta_k(2s,\varphi)$, удовлетворяющих уравнению
$B^{2s}(\zeta)=e^{i\varphi}$, имеем
$$
\|R\|_{L_2(\gamma_r)}^2=\int_{\gamma_r}|R(\zeta)|^2|d\zeta|=\frac{\pi
r}{s}\sum_{k=1}^{2s n}\frac{|R(\zeta_k)|^2}{\mu(\zeta_k)}.
 \eqno{(4)}
$$
При этом, в отличие от теоремы $1$, достаточно, чтобы сумма
кратностей полюсов в точках $z_1=0$ и $z_1^{*}=\infty$ не
превосходила $2s-1$. В частности, $(4)$ справедливо для
многочленов степени не выше $2s-1$ и функции $\mu$ $($тогда
$s=1$$)$:
$$
 \|\mu\|_{L_2(\gamma_r)}^2={\pi r}\sum_{k=1}^{2
n}{\mu(\zeta_k(2,\varphi))}.
 \eqno{(5)}
$$
}

\medskip
{\bf Примечание.} Из теоремы 2 получаются квадратурные формулы и
для $L_{2m}$-норм ($m=2,3,\ldots$)
$$
\|R\|_{L_{2m}(\gamma_r)}^{2m}=\int_{\gamma_r}|R(\zeta)|^{2m}|d\zeta|=\frac{\pi
r}{ms}\sum_{k=1}^{2\,m s
n}\frac{|R(\zeta_k)|^{2m}}{\mu(\zeta_k)},\qquad
\zeta_k=\zeta_k(2ms,\varphi),
 \eqno{(6)}
$$
применением (4) к рациональным функциям $R^{m}(z)$ с заменой $s$
на $ms$. Это же относится и к нижеследующим аналогичным формулам.

\medskip
{\bf 1.2.} Из теорем 1 и 2 в некоторых случаях можно получить
точные квадратурные формулы на множествах, которые являются
образами окружностей при отображениях рациональными функциями.
Сформулируем такой результат для отрезка $[-1,1]$. Пусть $W_n$
---  не пересекающееся с отрезком $[-1,1]$ конечное множество
(попарно различных точек) вида
$$
W_n=\{\infty,w_2,\ldots,w_{n}\}\cup
\{\overline{w_2},\ldots,\overline{w_{n}}\},\quad n\ge 2, \qquad
W_1=\{\infty\}
$$
(кратность вещественных точек объединения не учитывается). При
отображении Жуковского $w=(v+1/v)/2$ множеству $W_n$ соответствует
прообраз $V_{\nu}={\cal Z}_{\nu}\cup {\cal Z}_{\nu}^*$ вида (1),
содержащий $2\nu$ различных точек с $\nu\le 2n-1$ (по-прежнему
$z_1=0$, $z_1^{*}=\infty$). Положим (ср. с (2))
$$
B_0(z)=\prod_{k=1}^{\nu}\frac{z-z_k} {1-\overline{z_k} z},\qquad
\mu_0(\zeta)=\zeta\;\frac{B_0'(\zeta)}{B_0(\zeta)}>0,\qquad
|\zeta|=1.
$$

\medskip
{\bf Теорема 3.}  {\it Пусть $s\in \mathbb{N}$ и $R$
--- рациональная функция, все полюсы которой лежат на
$W_n$, причем отличные от $\infty$ полюсы имеют кратности не выше
$s$, а в бесконечно удаленной точке может быть полюс кратности не
выше $s-1$. Тогда при каждом $\varphi\in \mathbb{R}$
$$
\int_{-1}^{1}\frac{R(x)}{\sqrt{1-x^2}}\,dx=\frac{\pi}{s}
\sum_{k=1}^{s\nu} \frac{R(x_k(\varphi))}{\mu_0(\zeta_k(\varphi))},
\qquad x_k=\frac{1}{2}\left(\zeta_k+\frac{1}{\zeta_k}\right),
 \eqno{(7)}
$$
$$
\|R\|_{*,2}^2:=\int_{-1}^{1}\frac{|R(x)|^2}{\sqrt{1-x^2}}\,dx=
\frac{\pi}{2s} \sum_{k=1}^{2s \nu}
\frac{|R(\tau_k(\varphi))|^2}{\mu_0(\xi_k(\varphi))}, \qquad
\tau_k=\frac{1}{2}\left(\xi_k+\frac{1}{\xi_k}\right),
 \eqno{(8)}
$$
где $\nu\le 2n-1$, а узлы $\zeta_k=\zeta_k(\varphi)$ и
$\xi_k=\xi_k(\varphi)$ определяются соответственно из уравнений
$B_0^{s}(\zeta)=e^{i\varphi}$ и $B_0^{2s}(\xi)=e^{i\varphi}$.}

\medskip
{\bf Примечание.} При $\nu=1$, $W_1=\{\infty\}$, $\mu_0\equiv1$,
$\varphi=\pi$ и при замене $s$ на $2s$ точки насечки в (7)
определяются из уравнения $\zeta^{2s}=-1$. В этом случае точки
$x_k\in (-1,1)$ составляют двукратное множество нулей многочлена
Чебышева $T_s(x)$, в частности, оно расположено симметрично
относительно $x=0$. Отсюда получается известная квадратурная
формула Гаусса-Чебышева (см., напр., \cite[Гл. 10]{Krylov},
\cite{OSIPENKO}), точная на многочленах степени $\le 2s-1$:
$$
\int_{-1}^{1}\frac{f(t)}{\sqrt{1-t^2}}\,dt\approx
\frac{\pi}{2s}\sum_{k=1}^{2s} f(x_k)=\frac{\pi}{s}\sum_{k=1}^s
f(\lambda_k), \qquad \lambda_k=\cos \frac{2k-1}{2 s}\pi.
$$

\medskip
Сформулируем результат для действительной оси ${\mathbb R}$. Пусть
множество ${\cal F}_n$ соcтоит из $2n$ различных точек, не лежащих
на ${\mathbb R}$, и имеет вид
$$
{\cal F}_n=\{z_1,z_2,\ldots,z_{n}\}\cup
\{\overline{z_1},\overline{z_2},\ldots,\overline{z_{n}}\},\qquad
z_k\in {\mathbb C}^{+},\qquad z_k\ne\infty.
$$
Положим
$$
B_1(z)=\prod_{k=1}^{n}\frac{z-z_k}{z-\overline z_k}, \quad
\mu_1(x)=\frac{1}{2
i}\frac{B'_1(x)}{B_1(x)}=\sum_{k=1}^{n}\frac{{\rm
Im}\,z_k}{|x-\overline z_k|^2},\quad x\in {\mathbb R}.
  \eqno{(9)}
$$

\medskip
{\bf Теорема 4.} {\it Пусть $s\in \mathbb{N}$. Тогда для любой
правильной рациональной дроби $R$, все полюсы которой лежат на
множестве ${\cal F}_n$ и имеют кратности не выше $s$, при любом
$\varphi\in (0,2\pi)$ имеем
$$
\|R\|_{L_2(\mathbb{R})}^2=\frac{\pi}{2 s}\sum_{k=1}^{2s
n}\frac{|R(x_k(\varphi))|^2}{\mu_1(x_k(\varphi))},\qquad
\|\mu_1\|_{L_2(\mathbb{R})}^2=\frac{\pi}{2}\sum_{k=1}^{2
n}\,{\mu_1(t_k(\varphi))},
 \eqno{(10)}
$$
где в первом равенстве $x_k(\varphi)$ --- $($вещественные$)$ точки
насечки, удовлетворяющие уравнению $B_1^{2s}(x)=e^{i\varphi}$, а
во втором точки $t_k(\varphi)$ удовлетворяют уравнению
$B_1^{2}(t)=e^{i\varphi}$. }

\medskip
Из квадратурных формул для рациональных функций $R$ получаются
точные (экстремальные) и точные по порядку (относительно степени
$\deg R$) неравенства разных метрик (типа неравенств
С.~М.~Ни\-коль\-ско\-го), см. \S 3. Например, из теоремы 2
получается

\medskip
{\bf Теорема 5.} {\it Пусть $s\in \mathbb{N}$ и рациональная
функция $R$ удовлетворяет условию теоремы 2 $($конечные полюсы
имеют кратность не выше $s$, а в точках $z_1=0$ и $z_1^{*}=\infty$
сумма кратностей полюсов не превосходит $2s-1$$)$. Тогда
$$
\frac{|R(\zeta)|^2}{\mu(\zeta)}\le \frac{s} {\pi r}
\|R\|_{L_2(\gamma_r)}^2, \qquad \zeta\in \gamma_r.
  \eqno{(11)}
$$
Неравенство $(11)$ является точным: для любого $\zeta_0\in
\gamma_r$ найдется рациональная функция $R(\zeta_0,\zeta)$, для
которой оно обращается в равенство при $\zeta=\zeta_0$. В
частности, существует рациональная функция $R^*$, для которой
$$
\|R^*\|^2_{L_{\infty}(\gamma_r)}=\frac{s} {\pi r}
\|R^*\|_{L_2(\gamma_r)}^2 \|\mu\|_{L_{\infty}(\gamma_r)}.
  \eqno{(12)}
$$
Кроме того,
$$
\|R\|_{L_{\infty}(\gamma_r)} \le \frac{s}{\pi
r}\|R\|_{L_2(\gamma_r)}\|\mu\|_{L_2(\gamma_r)}.
  \eqno{(13)}
$$
Множитель перед произведением норм в $(13)$ нельзя уменьшить. }

\medskip
{\bf 1.3.} {\bf Примечания и примеры.} Подобные (3), (7)
квадратурные формулы, основанные на различных интерполяционных
тождествах (Бернштейна, Чебышева–Мар\-ко\-ва и др.) на
вещественной оси и окружности, были получены многими авторами (см.
работы \cite{ROVBA1}-\cite{OSIPENKO} и библиографию в них).
Результаты, близкие по форме к теоремам 1-3, получены в работах
Е.А. Ровбы, К.А. Смотрицкого, Е.В. Дирвука, Г. Мина. Например, для
рациональной функции $R(z)$ с симметричным относительно $\mathbb
R$ множеством простых полюсов $z_k$, $k=\overline{1,n-1}$
(допускается и полюс не выше первого порядка в бесконечно
удаленной точке), имеем \cite{ROVBA2}
$$
\int_{-1}^{1}\frac{R(x)}{\sqrt{1-x^2}}\,dx=\pi
\widetilde\sum_{k=1}^{n+1} \frac{R(t_k)}{\lambda_n(t_k)},\qquad
\lambda_n(x):=\sum_{k=1}^{n}\frac{\sqrt{1-z_k^2}}{1+xz_k},
 \eqno{(14)}
$$
где знак $\tilde{}$ означает, что первое и последнее слагаемое
следует разделить на 2, действительная часть радикалов берется
положительной, а в качестве узлов $t_2,\ldots,t_{n}$ берутся нули
синус-дроби Чебышева–Маркова
$$
N_n(x)=\frac{\sin \nu_n(x)}{\sqrt{1-x^2}},\qquad \hbox{где} \qquad
\nu_n(x)=-\int_{0}^{x}\frac{\lambda_n(t)}{\sqrt{1-t^2}}\, dt,
$$
и дополняются точками $t_1=-1$, $t_{n+1}=1$ (формулы с двумя
фиксированными узлами $\pm 1$ называются квадратурными формулами
типа Лобатто).

Формула (7) дополняет формулу (14) в следующем смысле. Несмотря не
внешнее сходство, эти формулы различаются: множества узлов и весов
в них имеют разную природу. Формула (14) работает на рациональных
функциях, которые в точке $z=\infty$ могут иметь только простой
полюс \cite{ROVBA2}. Кроме того, для нас важно, что узлы в (7), в
отличие от (14), переменные, и каждый узел $\zeta_k(\varphi)$
непрерывно пробегает (дважды) отрезок $[-1,1]$ при изменении
$\varphi$ от $0$ до $2\pi s \nu$. Это существенно используется
далее в неравенствах разных метрик.

Примерно тем же отличаются формулы (3) и (7) от формул В.Н.Русака,
Н.В.Гриба, Н.К. Филипповой (\cite{RUSAK1}, \cite{RUSAK2}) на
отрезке, а также от формул, полученных ими в случае окружности для
тригонометрических рациональных функций.

Отметим еще, что в \cite{D2010} получены сходные с (10)
квадратурные формулы для $L_2$-норм так называемых наипростейших
дробей (НД) порядка $n$:
 $$
\rho_n(z)=\sum_{k=1}^{n} (z-z_k)^{-1}.
    \eqno{(15)}
$$
Именно, если все $z_k\in \mathbb{C}^+$, и $x_k=x_k(\varphi)$ ---
точки насечки, удовлетворяющие уравнению $B_1^{2}(x)=e^{i\varphi}$
($k=\overline{1,2n}$), то
$$
\|\rho_n\|^2_2
 =\pi\sum\nolimits_{k=1}^{2n}\frac{({\rm
Re}\,\rho_n(x_k))^2}{\mu_1(x_k)}
=\pi\sum\nolimits_{k=1}^{2n}{\mu_1(x_k)},\qquad \mu_1(x)={\rm
Im}\,\rho_n(x).
     \eqno{(16)}
$$

\medskip
{\bf Примеры.} Приведем два примера применения (8) и (10). Если
взять $R(x)=(z-i)^{-s}$, то при $s=\overline{1,6}$
$$
\|R\|_{*,2}^2=\frac{\pi}{3\sqrt{2}}
 \left(|R(1)|^2+
 \frac{8}{5}\left|R\left(\sqrt{{3}/{5}}\right)\right|^2+
 \frac{8}{7}\left|R\left(\sqrt{{1}/{7}}\right)\right|^2+
 \frac{1}{2}|R(0)|^2\right)
$$
Положим $a= 2+ \sqrt {3}$, $b=2-\sqrt {3}$. Для любой правильной
рациональной дроби $R$ с полюсами кратности не выше 3,
расположенными в точках $\pm i$, из (10) получим
$$
 \|R\|_{L_2(\mathbb{R})}^2=\frac{\pi}{3}
\Bigl(\left| R\left( 1 \right) \right|^{2}+\left|{R}\left(-1
\right) \right|^{2}+
 2a \left(\left| R\left(a\right)\right|^{2}
 +\left| R\left(-a\right)\right|^{2}\right)+
 $$
 $$
+ 2b  \left(\left| {R}
 \left(b\right)  \right|^{2}+
 \left| R\left(-b\right)  \right|
^{2} \right)\Bigr).
$$

\medskip
{\bf\large 2. Доказательство теорем 1-4}
\medskip

\medskip
{\bf Доказательство теоремы 1.}  Сначала проведем доказательство в
случае простых полюсов $z_k$, т.е. при $s=1$. В этом случае по
условию точки $0$ и $\infty$ не могут быть полюсами дроби $R$.
Пусть
$$
R(z)=\sum_{k=2}^{n}\left(\frac{\alpha_k}{z-z_k}
+\frac{\beta_k}{r^2-z\overline{z_k}}\right),\qquad
F(z):=\frac{R(z)}{ B(z)-e^{i\varphi}}.
 $$
По теореме Коши о вычетах имеем
$$
\int_{|\zeta|=r}R(\zeta)|d\zeta|
=-ir\int_{|\zeta|=r}R(\zeta)\frac{d\zeta}{\zeta}
=-ir\int_{|\zeta|=r}
\sum_{k=2}^{n}\frac{\beta_k}{r^2-\zeta\overline{z_k}}
\frac{d\zeta}{\zeta} =\frac{2\pi}{ r} \sum_{k=2}^{n}\beta_k.
 \eqno{(17)}
$$
Рассмотрим разложение на простейшие дроби функции $F(z)$. Она
имеет простые полюсы в точках насечки $\zeta_k$ $(k=1,\ldots,n)$ и
в точках ${\cal Z}_n$ и устранимые особенности в точках ${\cal
Z}^{*}_n$. Поэтому разложение с учетом (2) и равенства
$B(\zeta_k)=e^{i\varphi}$ имеет вид
$$
F(z)=\sum_{k=1}^{n}\frac{R(\zeta_k)}{B'(\zeta_k)(z-\zeta_k)}
-\sum_{k=2}^{n}\frac{\alpha_k e^{-i\varphi}}{z-z_k} =
\sum_{k=1}^{n}\frac{{e^{-i\varphi}}\zeta_k
{R(\zeta_k)}}{\mu(\zeta_k)(z-\zeta_k)}
-\sum_{k=2}^{n}\frac{\alpha_k e^{-i\varphi}}{z-z_k}.
  \eqno{(18)}
$$
После подстановки $z=0$ в (18) с учетом $B(0)=0$, $R(0)\ne\infty$,
находим
$$
F(0)={R(0)}{e^{-i\varphi}}={e^{-i\varphi}}
\sum_{k=1}^{n}\frac{R(\zeta_k)}{\mu(\zeta_k)}
-e^{-i\varphi}\sum_{k=2}^{n}\frac{\alpha_k
}{z_k}\qquad\Rightarrow\qquad \sum_{k=2}^{n}\frac{\beta_k}{r^2} =
\sum_{k=1}^{n}\frac{R(\zeta_k)}{\mu(\zeta_k)}.
$$
Отсюда и из (17) получаем (3) при $s=1$.

К разобранному случаю сводится общий случай кратных полюсов,
отличных от $z_1=0$ и $z_1^{*}=\infty$. Действительно, пусть
$R(z)$ --- рациональная функция, полюсы которой лежат на ${\cal
Z}_n\cup {\cal Z}_n^*$, отличны от $z_1$ и $z_1^{*}$ и имеют
кратности не выше $s$. При $\nu=1,2,\ldots$ рассмотрим
последовательность множеств ${\cal Z}(\nu):=\{z_{\nu;k,j}\}$
$(j=\overline{1,s}, k=\overline{1,n})$ из попарно различных точек
со свойствами $\lim_{\nu\to\infty}z_{\nu;k,j}\to z_k$,
$z_{\nu;1,1}=z_1=0$. Тогда для функций $B(\nu;z)$ и $\mu(\nu;z)$
вида (2), построенных по множествам ${\cal Z}(\nu)\cup {\cal
Z}^{*}(\nu)$, при $\nu\to\infty$ имеем $B(\nu;\zeta)\to
B^s(\zeta)$ и $\mu(\nu;\zeta)\to s\mu(\zeta)$ равномерно на
$\gamma_r$. Кроме того, найдется последовательность правильных
рациональных функций $R(\nu;\zeta)$ с простыми отличными от $z_1$
и $z_1^{*}$ полюсами, лежашими на ${\cal Z}(\nu)\cup {\cal
Z}^{*}(\nu)$, для которой $R(\nu;\zeta)\to R(\zeta)$ равномерно на
окружности $\gamma_r$. Поэтому, применив уже доказанное к
$R(\nu;\zeta)$, $B(\nu;\zeta)$, $\mu(\nu;\zeta)$, затем предельным
переходом при $\nu\to\infty$ получим (3) для случая кратных
отличных от $z_1=0$ и $z_1^{*}=\infty$ полюсов.

Для завершения доказательства (3) и доказательства (4) понадобится

\medskip
{\bf Лемма 1.} {\it В условиях теоремы $1$ при
$\zeta_k=\zeta_k(s,\varphi)$ справедливо равенство }
$$
S(Q):=\frac{1}{s}\sum_{k=1}^{s
n}\frac{Q(\zeta_k)}{\mu(\zeta_k)}=q_0,\qquad \hbox{где} \qquad
Q(z):=\sum_{j=1-s}^{s-1}q_j z^j,\quad q_j\in {\mathbb C}.
 \eqno{(19)}
$$
Действительно, при $\sigma_j(z)=z^{j}$ и отрицательных $j\ge 1-s$
по аналогии с (18) находим
$$
\frac{\sigma_j(z)}{
B^s(z)-e^{i\varphi}}+\frac{\sigma_j(z)}{e^{i\varphi}}
=\frac{{e^{-i\varphi}}}{s}\sum_{k=1}^{n\,s}\frac{\zeta_k
{\sigma_j(\zeta_k)}}{\mu(\zeta_k)(z-\zeta_k)},
 $$
где предел при $z\to 0$ от левой части равен нулю (т.к. $z_1=0$
является нулем порядка $s>|j|$ для $B^s$). Отсюда при $z\to 0$
получаем, что сумма $S(\sigma_j)$ равна нулю. Для положительных
$j$ также $\overline{S(\sigma_j)}=r^{2j}S(\sigma_{-j})=0$.
Аналогично разложением дроби $(B^s(z)-e^{i\varphi})^{-1}$
получается равенство {$S(\sigma_0)=1$}.

Интеграл в (3) от степеней $\sigma_j(z)$, $1-s\le j\le s-1$, $j\ne
0$, равен нулю, поэтому из (19) следует, что добавление к $R$
линейной комбинации таких степеней не изменяет обе части равенства
(3). Для константы $\sigma_0=1$ обе части равны $2\pi r$. Этим
завершается доказательство (3). $\square$

\medskip
{\bf Доказательство теоремы 2.} Рассмотрим функцию
$R_0(z)=R(z)\overline{R(z)}$. При $z=\zeta\in \gamma_r$ имеем
$$
R_0(\zeta)=|R(\zeta)|^2=R(\zeta)\overline{R(\zeta)}=
R(\zeta)\overline{R}(\overline{\zeta})=
R(\zeta)\overline{R}\left(\frac{r^2}{\zeta}\right),
$$
где $\overline{R}(\cdot)$ означает операцию сопряжения
коэффициентов функции $R$. Считая теперь $\zeta$ свободной
переменной на ${\mathbb C}$, несложно видеть, что $R_0$ как
функция от $\zeta$ является рациональной функцией вида
$$
R_0(\zeta)=P_m(\zeta)+Q_{\nu}(1/\zeta)+R_1(\zeta),
$$
где $P_m$ и $Q_{\nu}$ --- многочлены, степени $m$ и ${\nu}$
которых не превосходят суммы кратностей полюсов функции $R$ в
точках $z_1=0$ и $z_1^{*}=\infty$, т.е. $m\le 2s-1$, ${\nu}\le
2s-1$, $R_1$ --- правильная рациональная дробь, ненулевые полюсы
которой лежат на множестве ${\cal Z}_n\cup {\cal Z}_n^*\setminus
\{\infty\}$ и имеют кратность не выше $2s$. Применив теорему 1 с
заменой $s$ на $2s$, получим формулу (3) для $R_1$. Эта же формула
верна и для многочленов $P_m$ и $Q_{\nu}$. Это следует из леммы 2
(где $s$ заменяется на $2s$) и из теоремы о среднем:
$\int_{\gamma_r}P_m(z)|dz|=2\pi r P_m(0)$,
$\int_{\gamma_r}Q_{\nu}(z)|dz|=2\pi r Q_{\nu}(0)$. $\square$

\medskip
{\bf Доказательство теоремы 3.}  При  $x\in[-1,1]$ сделаем замену
Жуковского
$$
x=\frac{1}{2}\left({\zeta}+\frac{1}{{\zeta}}\right)=\cos t,\quad
{\zeta}=e^{it},\quad\sqrt{1-x^2}=|\sin{t}|,\quad d\,
x=-\sin{t}\;d\,t.
 \eqno{(20)}
$$
Учитывая, что при обходе точки ${\zeta}$ единичной окружности
точка $x$ пробегает отрезок $[-1,1]$ дважды, получим
$$
2\int_{-1}^{1}\frac{R(x)}{\sqrt{1-x^2}}\,dx= \int_{\gamma_1}
R_0({\zeta}) |d \zeta|,\qquad \hbox{где} \qquad
R_0(z):=R\left(\frac{1}{2}\left(z+\frac{1}{z}\right)\right).
  \eqno{(21)}
$$
Отличные от $z_1=0$ и $z_1^{*}=\infty$ полюсы рациональной дроби
$R_0(z)$ имеют кратность не выше $s$, а в точках $z_1=0$ и
$z_1^{*}=\infty$ могут быть полюсы кратности не выше $s-1$ .
Поэтому утверждение теоремы 3 следует из квадратурных формул (3) и
(4) с заменой $n$ на $\nu$:
$$
2\int_{-1}^{1}\frac{R(x)}{\sqrt{1-x^2}}\,dx=\frac{2\pi}{s}
\sum_{k=1}^{s \nu}\frac{R_0(\zeta_k)}{\mu(\zeta_k)},\quad
2\|R\|_{*,2}^2=\frac{\pi}{s} \sum_{k=1}^{2s\cdot
\nu}\frac{|R_0(\zeta_k)|^2}{\mu(\zeta_k)}.
$$

\medskip
{\bf Доказательство теоремы 4.} Для удобства доказательства будем
считать $z_1=i$. При замене $x=-i(\zeta+1)(\zeta-1)^{-1}$,
$x\in{\mathbb R}$, $\zeta\in \gamma_1$, получаем
$$
\|R\|_{L_2(\mathbb{R})}^2=2\int_{\gamma_1}
\left|R_0(\zeta)\right|^2 |d\zeta|,\quad
R_0(\zeta)=\frac{1}{\zeta-1}
R\left(i\frac{1+\zeta}{1-\zeta}\right),\quad
\zeta=\frac{x-i}{x+i},
    \eqno{(22)}
$$
где $R_0$ --- рациональная дробь, полюсы которой лежат на
симметричном относительно $\gamma_1$ множестве
$$
\{\xi_k\}\cup \{\xi^{*}_k\},\qquad \hbox{где} \qquad
\xi_k=\frac{z_k-i}{z_k+i}, \qquad k=1,\ldots,n,
$$
и во всех точках, отличных от $\xi^{*}_1=\infty$, имеют кратность
не выше $s$, а кратность полюса в точке $\xi^{*}_1=\infty$, как
видно из (22), не выше $s-1$. В точке $z=1$ --- устранимая
особенность (т.к. $R(\infty)=0$ по условию). Следовательно, $R_0$
удовлетворяет условиям теоремы 2. Далее, легко проверить, что при
указанной замене для произведений Бляшке (2) и (9) имеем связь
$B(\zeta)=A\cdot B_1(x)$, т.е.
$$
\prod_{k=1}^{n}\frac{\zeta-\xi_k}{1-\zeta\overline{\xi_k}}= A
\prod_{k=1}^{n} \frac {x-z_k}{x-\overline{z_k}},\qquad \hbox{где}
\qquad A:=\prod_{k=1}^{n} \frac {i-\overline{z_k}} { i+z_k}.
$$
Отсюда с учетом (2) и (9) находим
$$
\mu(\zeta)=\zeta\frac{B'(\zeta)}{B(\zeta)}=
\zeta\frac{B_1'(x)}{B(x)}\frac{2i}{(1-\zeta)^2}=
-\frac{4\zeta}{(1-\zeta)^2}\,\mu_1(x),
$$
Наконец, учитывая (22) и (4), получаем (10):
$$
2\|R_0\|_{L_2(\gamma_1)}^2
 =\frac{2\pi
r}{s}\sum_{k=1}^{2s n}\frac{|R_0(\zeta_k)|^2}{\mu(\zeta_k)}
=\frac{\pi}{2s}\sum_{k=1}^{2s n}\frac{|R(x_k)|^2}{
\mu_1(x_k)},\quad x_k=-i\frac{\zeta_k+1}{\zeta_k-1}.
$$

\bigskip
\begin{center}
{\bf\large 3. Неравенства разных метрик для рациональных функций
\\ и многочленов}
\end{center}

\medskip
Неравенства разных метрик для алгебраических и тригонометрических
многочленов хорошо известны. Они восходят к работам Д.Джексона
\cite{Jackson} и С.М.\,Ни\-коль\-ско\-го \cite{Nikolskiy}.
Дальнейшее развитие эта тематика получила благодаря работам
Н.К.\,Бари, С.Б. Стечкина, В.В. Арестова, Л.В.Тайкова, Г.\,Сегё,
А.\,Зигмунда, П.Л.\,Ульянова, М.К.\,По\-та\-по\-ва,
А.Ф.\,Ти\-ма\-на, И.И.\,Иб\-ра\-ги\-мова, Дж.И. Мамедханова  и
многих других авторов (все источники здесь невозможно перечислить,
мы ограничимся ссылкой на работы \cite{Nikolskiy1}, \cite{Bari},
\cite{Arestov}, \cite{KONIAGIN}, \cite{Mamedhanov}, \cite{Timan} и
на соответствующую библиографию в них; см. также ссылки ниже). В
качестве примера приведем неравенство для алгебраических
многочленов $P_n$ степени $\le n$ на ограниченном промежутке
$E\subset {\mathbb R}$:
$$
 \|P_n\|_{L_p(E)}\le C(E,p,r)
n^{2\left(\frac{1}{r}-\frac{1}{p}\right)} \|P_n\|_{L_r(E)},\qquad
\infty\ge p>r\ge 1.
    \eqno{(23)}
$$
Сходные неравенства верны и для тригонометрических многочленов,
причем в случае $E=[-\pi,\pi]$ множитель 2 в показателе степени
можно заменить единицей (см. \cite{Bari},\cite{Nikolskiy}).

Некоторые неравенства разных метрик для НД (15) были получены в
\cite{DAN-DOD}, \cite{[DAN1]}. Приведем одно такое неравенство:
$$
 \left\| \rho_n \right\|_{L_p(\mathbb R)}^{q} \leq A(p,r)
\left\| \rho_n \right\|_{L_r(\mathbb R)}^{s},\qquad
\hbox{где}\qquad p^{-1}+q^{-1}=1,\quad r^{-1}+s^{-1}=1.
  \eqno{(24)}
$$
Отметим, что, в отличие от (23), неравенства (24) для НД нелинейны
относительно сравниваемых норм. Вообще, аналогичная нелинейность
возникает и в других неравенствах для НД, кроме того неравенства
для НД содержательны и на бесконечных промежутках и при
произвольном соотношении между $p$ и $r$ (при $p<r$ они также не
являются тривиальными и не следуют из неравенства Гельдера даже на
конечных промежутках) \cite{DAN-DOD}.

В работах В.Ю.\,Прота\-со\-ва, И.Р.\,Каю\-мо\-ва, А.В.\,Каюмовой,
В.И.\,Дан\-чен\-ко \cite{Protasov}--\cite{Kayumova} были получены
оценки другого типа. В этих работах $L_p$-нормы для НД в основном
оценивались посредством специальных сумм, явно зависящих от
полюсов $z_k$. Например, в \cite{Protasov} показано, что
$$
\|\rho_n\|_{{L_p(\mathbb R)}}\ge
\gamma_p\left(\sum_{k=1}^n\frac{1}{|{\rm
Im}\,z_k|^{p-1}}\right)^{1/p}, \qquad p>1,
$$
с определенным точным множителем $\gamma_p$, зависящим лишь от
$p$. В \cite{D2010} показано, что если все $z_k\in
\mathbb{C}^{+}$, то при $1<p<3$ имеем двустроннюю оценку
$$ \|\rho_n\|_{L_p(\mathbb R)}^p\;
\cos\frac{\pi(p-2)}{2}\le 2\pi\;{\rm
Im}\,\left(e^{-i\frac{\pi(p-2)}{2}}\sum_{k=1}^{n}
\rho_n^{p-1}(\overline{z_k})\right)\le \|\rho_n\|_{L_p(\mathbb
R)}^p.
$$
Возникает интерес к распространению неравенств типа (24) на другие
метрики и множества. Некоторые результаты в этом направлении
получены в \cite{[DAN-DAN]}. В следующем параграфе получены точные
оценки $\sup$-норм рациональных функций и, в частности,
многочленов и напростейших дробей через их интегральные нормы в
случае окружностей, прямых и их сегментов.

\medskip
{\bf 3.1. Случай окружности.} Докажем теорему 5, а затем как ее
частный случай получим некоторые неравенства для многочленов.

\medskip
{\bf 3.1.1. Доказательство теоремы 5. } Выбрав в (4) параметр
$\varphi$ так, чтобы в точке насечки $\zeta_1=\zeta_1(2s,\varphi)$
функция $|R(\zeta)|^2/\mu(\zeta)$ принимала наибольшее значение,
получим (11):
$$
\max_{\zeta\in
\gamma_r}\frac{|R(\zeta)|^2}{\mu(\zeta)}=\frac{|R(\zeta_1)|^2}{\mu(\zeta_1)}\le
\frac{s} {\pi r} \|R\|_{L_2(\gamma_r)}^2
$$
Далее, из неравенства Коши-Бу\-ня\-ков\-ского и (4) имеем
$$
\sum_{k=1}^{2s n}|R(\zeta_k)| \le \sqrt{\sum_{k=1}^{2s
n}\frac{|R(\zeta_k)|^2}{\mu(\zeta_k)}}\sqrt{\sum_{k=1}^{2s
n}{\mu(\zeta_k)}}=\frac{s}{\pi
r}\|R\|_{L_2(\gamma_r)}\|\mu\|_{L_2(\gamma_r)},
 $$
где $\zeta_k=\zeta_k(2s,\varphi)$, откуда, выбрав в (4) параметр
$\varphi$ так, чтобы в точке $\zeta_1=\zeta_1(2s,\varphi)$ функция
$|R(\zeta)|$ принимала наибольшее значение, получим (13).

Конечно, выбор в сумме (4) лишь одного слагаемого может
значительно огрублять оценки, особенно когда слагаемые
незначительно отличаются друг от друга. Тем не менее оценка (11)
является точной на классе допустимых рациональных функций.
Покажем, например, что равенство в (11) достигается на
рациональных функциях вида
$$
R(z)=R(s,n,\varphi;z)=\frac{B^{s}(z)-e^{i
\varphi}B^{-s}(z)}{z-\zeta_1},\qquad \zeta_1=\zeta_1(2s,\varphi),
   \eqno{(25)}
$$
где $B$ --- произведение Бляшке из (2). Очевидно, функция $R$
удовлетворяет условиям теоремы 2 и, следовательно, для нее
справедливо равенство (4). С другой стороны, во всех точках
насечки, отличных от $\zeta_1$, имеем $R(\zeta_k)=0$, $k\ne 1$.
Это означает, что функция $|R(\zeta)|^2/\mu(\zeta)$ принимает
наибольшее значение в точке $\zeta_1$ и в (11) получается
равенство. Для доказательтва (12) достаточно дополнительно
потребовать, чтобы при соответствующем $\varphi$ в точке
$\zeta_1=\zeta_1(2s,\varphi)$ функция $\mu$ принимала наибольшее
значение. $\square$

В дополнение к теореме 5 заметим, что норма $\mu$ в (13) легко
подсчитывается:
$$
\|\mu\|^2_{L_2(\gamma_r)}=r^2\int_{\gamma_r}|B'(\zeta)|^2 |d
\zeta|=-i r^2\int_{\gamma_r}|B'(\zeta)|\frac{B'(\zeta)}{B(\zeta)}
d \zeta=-i
r\int_{\gamma_r}\zeta\left(\frac{B'(\zeta)}{B(\zeta)}\right)^2 d
\zeta=
$$
$$
=2\pi r\sum_{k=1}^{n}{\rm Res}_{z_k}\,
z\left(\frac{B'(z)}{B(z)}\right)^2=2\pi
r\left(n^2+2\sum_{j=2}^{n}\sum_{k=2}^{n} \frac{z_k \overline{z_j}
}{r^2-z_k \overline{z_j}}\right).
  \eqno{(26)}
$$
Отсюда следует, например, что множитель перед произведением норм в
(13) нельзя уменьшить. Действительно, при $R\equiv \mu$ и $s=1$ из
(13) и (26) имеем (см. также (5))
$$
\|\mu\|_{L_{\infty}(\gamma_r)}\le \frac{1}{\pi
r}\|\mu\|^2_{L_2(\gamma_r)}= 2n^2+4\sum_{j=2}^{n}\sum_{k=2}^{n}
\frac{z_k \overline{z_j} }{r^2-z_k \overline{z_j}},
  \eqno{(27)}
$$
а здесь, как легко проверить, при $z_2\nearrow r$ и прочих
фиксированных полюсах отношение левой и правой частей стремится к
единице. Отметим еще, что если все полюсы $R$ лежат во внешности
кольца $K_{\varepsilon}=\{z:\;\varepsilon r< |z|<
r/\varepsilon\}$, $\varepsilon\in (0,1)$, то, оценив $\mu(\zeta)$
сверху величиной $n(1+\varepsilon)/(1-\varepsilon)$, из (11)
получим
$$
\|R\|_{L_{\infty}(\gamma_r)} \le \sqrt{\frac{n s}{\pi
r}\,\frac{1+\varepsilon}{1-\varepsilon}}\;\|R\|_{L_2(\gamma_r)},
\qquad z_k\notin K_{\varepsilon}.
$$

\medskip
{\bf 3.1.2. Неравенства для многочленов.} Квадратурный метод
позволяет единообразно получить ряд неравенств для
тригонометрических и алгебраических многочленов. Приведем
некоторые оценки.

Пусть $T_{s-1}(t)$ --- тригонометрический многочлен степени не
выше $s-1$ с комплексными коэффициентами. Пусть $p>0$ и $m=m(p)$
-- наименьшее натуральное число, удовлетворяющее неравенству
$2m\ge p$. Тогда
$$
\|T_{s-1}\|_{L_{\infty}[0,2\pi]}\le \alpha(p)\cdot
s^{\frac{1}{p}}\,\|T_{s-1}\|_{L_{p}[0,2\pi]},\qquad
\alpha(p)\le\left(\frac{m}{\pi}\right)^{1/p}<\frac{11}{10},\qquad
p>0.
  \eqno{(28)}
$$
Аналогичное неравенство справедливо для алгебраического многочлена
$P_{2s-1}$ с комплексными коэффициентами степени не выше $2s-1$:
$$
\|P_{2s-1}\|_{L_{\infty}(\gamma_r)}\le \alpha(r,p)\cdot
s^{\frac{1}{p}}\, \|P_{2s-1}\|_{L_p(\gamma_r)},\qquad
\alpha(r,p)\le \left(\frac{m}{\pi r}\right)^{1/p},\qquad p>0.
 \eqno{(29)}
$$
Для доказательства (28) рассмотрим рациональную функцию
$R(z)=\sum_{k=1-s}^{s-1} \sigma_k z^k$ такую, что
$T_{s-1}(t):=R(e^{it})$.  В этом случае функция $R$ имеет два
полюса $z_1=0$ и $z_1^{*}=\infty$, сумма кратностей которых не
превосходит $2s-2$. Поэтому из примечания к теореме 2 при $n=1$,
$r=1$ и $\mu(e^{it})\equiv 1$ получается:
$$
\sum_{k=1}^{2s}|T_{s-1}(t_k)|^{2m}= \frac{s
m}{\pi}\|T_{s-1}\|^{2m}_{L_{{2m}}[0,2\pi]},
$$
где $t_k=t_k(\varphi)=\arg \zeta_k(2s m,\varphi)$. Отсюда, выбрав
$\varphi$ так, чтобы в точке $t_1$ функция $|T_{s-1}|$ принимала
наибольшее значение, получим
$$
\|T_{s-1}\|^{2m}_{L_{\infty}[0,2\pi]}=|T_{s-1}(t_1)|^{2m}\le
\frac{{s m}}{\pi}\|T_{s-1}\|^{2m}_{L_{{2m}}[0,2\pi]}\le \frac{s
m}{\pi}\|T_{s-1}\|^p_{L_{p}[0,2\pi]}
\|T_{s-1}\|^{{2m}-p}_{L_{\infty}[0,2\pi]},
$$
откуда и следует оценка (28). Неравенство (29) получается
аналогично из (6) при $R(z)=P_{2s-1}(z)$, $n=1$,
$\mu(e^{it})\equiv 1$.

\medskip
{\bf Замечания о точности.} Первое неравенство типа (28) для
вещественных тригонометрических многочленов получено
Д.Джек\-со\-ном и С.М. Никольским \cite{Jackson},\cite{Nikolskiy}
(при $p\ge 1$ c константой $\alpha(p)=2$). При $p=1$ поведение
наилучшего множителя $\alpha^{*}(1)$ в (28) для вещественных
многочленов изучали С.Б.Стечкин и Л.В.Тайков. Ими показано, что
$\alpha^{*}(1)=c+ o(1)$ при $s\to\infty$ с некоторым $0.539<c<
0.58$, см. \cite{Taikov2}.

Заметим, что на самом деле неравенство (28) справедливо в классе
${\cal T}_s$ многочленов $T_s$ степени $\le s$, которые
представляются в виде $\sum_{k=1-s}^{s} \sigma_k z^k$  или
$\sum_{k=-s}^{s-1} \sigma_k z^k$, $z=e^{it}$, $\sigma_k\in
{\mathbb C}$ (т.к. сумма кратностей полюсов в обоих случаях не
превосходит $2s-1$). В классе ${\cal T}_s$ неравенство (28) точно
при $p=2$ ($m=1$); например, легко проверяется, что при
$\alpha(2)= {\pi}^{-1/2}$ оно достигается на многочленах
$$
T^{*}_s(t)=1+2\sum_{k=1}^{s-1}\cos kt  +e^{s t
i}=1+\sum_{k=1}^{s-1}(z^k+z^{-k})+z^{s}, \quad z=e^{it}.
$$
В этом случае и левая, и правая части в (28) равны
$T^{*}_s(0)=2s$.

При $p=2$ ($m=1$) точным является и неравенство (29), оно
достигается на $P^{*}_{2s-1}(z)=(z^{2s}-r^{2s})(z-r)^{-1}$. В
остальных случаях вопрос о точности констант $\alpha(p)$ остается
открытым.

\medskip
Приведем еще одно неравенство, получающееся из (3) при $n=1$,
$r=1$ и $\mu(e^{it})\equiv 1$:
$$
\max_{\varphi}\left|\sum_{k=1}^{s}T_{s-1}(t_k(\varphi))\right|\le
\frac{s}{2\pi}\|T_{s-1}\|_{L_{1}[0,2\pi]},\qquad
t_k=t_k(\varphi)=\arg \zeta_k(s,\varphi),
  \eqno{(30)}
$$
откуда для {\it неотрицательных} тригонометрических многочленов
получается оценка вида (28) при $p=1$ с меньшим множителем
$\alpha(1)=1/(2\pi)$.

\medskip
{\bf 3.2. Случай отрезка.} Пусть $R(z)=P_{s-1}(z)$
--- алгебраический многочлен с комплексными
коэффициентами степени не выше $s-1$, $p>0$ и $m=m(p)$ ---
наименьшее натуральное число, удовлетворяющее неравенству $2m\ge
p$. Тогда аналогично (28) и (29) из (8) (при $n=1$,
$\mu_0(e^{it})\equiv 1$) с учетом примечания к теореме 2
получаем\footnote{Неравенства такого типа хорошо известны, см.,
например, \cite{POLYA}, \cite{Labelle}, \cite{Arestov-D}. Однако
авторам не известны такие неравенства для комплексных многочленов
с более точными константами.}
$$
\|P_{s-1}\|_{L_{\infty}([-1,1])}\le\beta(p)
s^{1/p}\left(\int_{-1}^{1}
\frac{|P_{s-1}(x)|^p}{\sqrt{1-x^2}}\,dx\right)^{1/p},\quad
\beta(p)\le\left(\frac{2m}{\pi}\right)^{1/p}<\frac{6}{5}.
 \eqno{(31)}
$$
Это неравенство является точным по порядку $s$ \cite{Arestov-D}.
Более того, при $p=2$ существуют вещественные многочлены
$P^{*}_{s-1}(x)$, для которых
$$
 |P^{*}_{s-1}(1)|= \frac{\sqrt{s}}{\sqrt{\pi
}}\|P^{*}_{s-1}\|_{*,2}.
   \eqno{(32)}
$$
Такое равенство выполняется, например, для многочленов
$P^{*}_{s-1}(x)$, которые при замене (20) принимают вид
$$
P^{*}_{s-1}(x)=C\frac{1}{\zeta^{s-1}}\frac{\zeta^{2s}-1}{\zeta^2-1},\qquad
x=\frac{1}{2}\left({\zeta}+\frac{1}{{\zeta}}\right).
  \eqno{(33)}
$$
Такие многочлены существуют, это многочлены Якоби
$P^{*}_{s-1}(x)=P^{(1/2,1/2)}(x)$ (см., напр., \cite{SEGE}),
удовлетворяющие дифференциальным уравнениям
$$
(x^2-1)y''(x)+3x y'(x)+(1-s^2)y(x)=0.
   \eqno{(34)}
$$
Действительно, при замене (33) для функции
$Y(\zeta)(\zeta^2-1)^{-1}=y(x)$ уравнение (34), как несложно
проверить, преобразуется к уравнению Эйлера
$$
\zeta^2\,Y''(\zeta)-\zeta\, Y'(\zeta)+(1-s^2)\, Y(\zeta)=0,
$$
общее решение которого имеет вид $\zeta^{1-s}
(c_1\zeta^{2s}+c_2)$. Значит, для полиномиального решения имеем
$c_2=-c_1$ и (33) выполнено. В этом случае применим (8) с
$\varphi=0$. Точки насечки $\xi_k$ определяются как корни
уравнения $\zeta^{2s}=1$ и поэтому в сумме (8) остается только два
равных ненулевых слагаемых $|P^{*}_{s-1}(\pm 1)|^2$. Отсюда и
следует (32).

\medskip
{\bf 3.3. Случай действительной оси.} Аналогично теореме 5 из (10)
получается

\medskip
{\bf Теорема 6.} {\it Пусть $s\in \mathbb{N}$ и правильная
рациональная дробь $R$ удовлетворяет условию теоремы $4$ $($полюсы
не лежат на $\mathbb R$ и имеют кратности не выше $s$$)$. Тогда
$$
\frac{|R(x)|^2}{\mu_1(x)}\le \frac{2s} {\pi }
\|R\|_{L_2(\mathbb{R})}^2,\qquad x\in \mathbb{R}.
$$
Это неравенство является точным: для любого $x_0\in \mathbb{R}$
найдется рациональная функция $R(x_0,x)$, для которой оно
обращается в равенство при $x=x_0$. В частности, существует
рациональная функция $R^*$, для которой
$$
\|R^*\|^2_{L_{\infty}(\mathbb{R})}=\frac{2s} {\pi}
\|R^*\|_{L_2(\mathbb{R})}^2 \|\mu_1\|_{L_{\infty}(\mathbb{R})}.
$$
Кроме того,
$$
\|R\|_{L_{\infty}(\mathbb{R})} \le \frac{2s}{\pi
}\|R\|_{L_2(\mathbb{R})}\,\|\mu_1\|_{L_2(\mathbb{R})},\qquad
\|\mu_1\|_{L_{\infty}(\mathbb{R})}\le \frac{2}{\pi
}\|\mu_1\|^2_{L_2(\mathbb{R})},
  \eqno{(35)}
$$
В $(35)$ множитель $2/\pi$ уменьшить нельзя, это проверяется на
рациональных функциях $R(z)\equiv\mu_1(z)$ с парой простых
полюсов.}

Экстремальные рациональные функции $R^*$ в теореме 6 строятся как
и в теореме 5 (см. (25)). Заметим, что $\mu_1(x)={\rm
Im}\,\rho_n(x)$ для НД (15) с $z_k=x_k+iy_k\in \mathbb{C}^+$, так
что если взять $R=\rho_n$ и $s=1$, то из (35) получается
$$
\|\rho_n\|_{L_{\infty}(\mathbb{R})} \le \frac{2}{\pi
}\|\rho_n\|_{L_2(\mathbb{R})}\,\|{\rm
Im}\,\rho_n\|_{L_2(\mathbb{R})},
$$
что уточняет сходный результат из \cite{DAN-DOD}, где вместо ${\rm
Im}\,\rho_n$ было $\rho_n$. Отметим еще, что по теореме Коши о
вычетах имеем (см., например, \cite{D2010})
$$
2^{-1}\|\rho_n\|^2_2=\|\mu_1\|^2_2=\pi
\sum\nolimits_{k,m=1}^{n}\frac{y_k+y_m}{(x_k-x_m)^2+(y_k+y_m)^2}.
$$

\medskip
{\bf 3.4. Неравенство разных метрик для НД.} Пусть ${R}(z)$ ---
рациональная функция, все полюсы которой лежат во внешности
окружности $\gamma_r$ и ${R}(\infty)=0$. Для НД (15) с полюсами,
лежащими внутри $\gamma_r$, положим
$$
\rho(z)=\rho_n(z)+{R}(z).
  \eqno{(36)}
$$
Докажем вспомогательное предложение о разделении особенностей. Для
краткости всюду ниже в этом параграфе применяем обозначение
$\|\cdot\|_{p}:=\|\cdot\|_{L_p(\gamma_r)}$.

\medskip
{\bf Лемма 2.} {\it На окружности $\gamma_r$, $r>0$, независимо от
порядка $n\ge 1$ имеем неравенство}
$$
\|\rho_n\|_{\infty}\le 3 {\eta}\ln (r {\eta}+1),\qquad \hbox{где}
\qquad {\eta}:= \|\rho\|_{\infty}.
 \eqno{(37)}
$$

\medskip
{\bf Доказательство.} Достаточно установить (37) в случае $r=1$;
общий случай получается заменой $\rho(z)$ на $r\,\rho(r\,z)$.
Имеем
 $$
 1\le n\le\frac{1}{2\pi}\int_{\gamma_1}
 |\rho(\zeta)||d\,\zeta|\le {\eta},\quad
 |\rho'_n(z)|\le
 \frac{{\eta}}{2\pi}\int_{\gamma_1}\frac{|d\zeta|}{|\zeta-z|^2}
 =\frac{{\eta}}{|z|^2-1}
  \eqno{(38)}
$$
при $|z|>1$. Будем для удобства считать
$z_1=1-\varepsilon\in(0,1)$. Пусть $z=\theta>1$. Из предыдущего
неравенства имеем
$$
|\rho_n(z)|\le {\eta}\int_{|z|}^{\infty}\frac{d
t}{t^2-1}=\frac{{\eta}}{2}\ln\frac{|z|+1}{|z|-1}\quad
\Rightarrow\quad \frac{1}{{\theta}-z_1}
\le
 \frac{{\eta}}{2}\ln\frac{{\theta}+1}{{\theta}-1}.
 \eqno{(39)}
$$
Итак, разность
$$
v({\eta},\varepsilon;{\theta})=
 \frac{{\eta}}{2}\ln\frac{{\theta}+1}{{\theta}-1}
 -\frac{1}{{\theta}-1+\varepsilon}
$$
должна быть положительной при всех ${\theta}>1$. Отсюда следует,
что $\varepsilon>(2{\eta}^2)^{-1}$. Действительно, при
$\varepsilon=(2{\eta}^2)^{-1}$, подставив ${\theta}=1+{\eta}^{-3}$
в $v({\eta},\varepsilon;{\theta})$, получим, как легко проверить,
отрицательное значение при всех ${\eta}\ge 1$ (из (38) следует,
что всегда ${\eta}\ge n\ge 1$). Итак, $1-|z_1|\ge
(2{\eta}^2)^{-1}$. То же верно для любого полюса $z_k$. Отсюда и
из (39) при любом $z$, $|z|>1$, и $\zeta=z/|z|$ получаем
$$
|\rho_n(\zeta)|\le |\rho_n(z)|+|\rho_n(z)-\rho_n(\zeta)|\le
\frac{{\eta}}{2}\ln\frac{|z|+1}{|z|-1}+
|z-\zeta|\sum_{k=1}^n\frac{1}{|z-z_k||\zeta-z_k|}
$$
$$
\le \frac{{\eta}}{2}\ln\frac{|z|+1}{|z|-1}+ |z-\zeta|\, {n}\cdot
2{\eta}^2\cdot 2{\eta}^2\le
\frac{{\eta}}{2}\ln\frac{|z|+1}{|z|-1}+ 4{\eta}^5\;(|z|-1).
$$
Если взять $|z|={2^{-1}}{{\eta}^{-2}}\,{\sqrt {1+4\,{\eta}^4}}$,
то отсюда получим нужное неравенство (37). $\square$

\medskip
Пусть ${\cal Z}=\{z_1,\ldots,z_{n}\}\cup \{0\}$, где $z_k$ полюсы
НД (15) (лежат внутри $\gamma_r$). Для удобства доказательства
следующей теоремы будем считать, что точки $z_k$ попарно различны
(общий случай сводится к этому предельным переходом <<слияния
полюсов>>, например, как это делалось при доказательстве теоремы
1). Для числа $\nu$ точек множества ${\cal Z}$ имеем $\nu\le n+1$.
Пусть ${R}(z)$ --- рациональная функция с ${R}(\infty)=0$, все
полюсы которой лежат на ${\cal Z}^*$ и имеют кратность не выше
$s$. Тогда справедлива

\medskip
{\bf Теорема 7.} {\it Для рациональной функции $(36)$ справедливы
неравенства $($обоз\-на\-чение: ${\eta}= \|\rho\|_{\infty})$
$$
\frac{{\eta}^2}{6r
 {\eta}\ln(1+{\eta})-(n-1)}
 \le
 \frac{s}{\pi
r}\|\rho\|_{2}^2,\qquad \frac{{\eta}}{\ln(1+{\eta})}
  \le
 \frac{6s}{\pi
}\|\rho\|_{2}^2.
  \eqno{(40)}
$$
В частности, при ${R}(z)\equiv 0$ имеем }
$$
\max_{|\zeta|=r}\frac{|\rho_n(\zeta)|^2}{2r
 |{\rho_n}(\zeta)|-(n-1)}\le
 \frac{1}{\pi
r}\|\rho_n\|_{2}^2,\qquad \|\rho_n\|_{\infty}\le
 \frac{2}{\pi }\|\rho_n\|_{2}^2.
 \eqno{(41)}
$$

\medskip
{\bf Доказательство.} По множеству ${\cal Z}$ определим как в (2)
функции $B$ и $\mu$. При $|\zeta|=r$ с учетом (2) и (37) имеем
$$
\mu(\zeta)\le 1+{\rm Re}\,(2\zeta\, {\rho}_{n}(\zeta))-n,\qquad
\mu(\zeta)\le 6r\, {\eta}\,\ln(1+{\eta})-n+1,
   \eqno{(42)}
 $$
откуда, учитывая (4), при соответствующем выборе $\varphi$ (как в
п. 3.1.1) получаем
$$
\frac{{\eta}^2}{6r
 {\eta}\ln(1+{\eta})-n+1}
  \le\sum_{k=1}^{2 n s}\frac{|\rho(\zeta_k)|^2}{\mu(\zeta_k)}=
 \frac{s}{\pi
r}\|\rho\|_{2}^2,\quad \frac{{\eta}}{\ln(1+{\eta})}
  \le
 \frac{6s}{\pi
}\|\rho\|_{2}^2.
$$
Если ${R}(z)\equiv 0$, то можно воспользоваться (4) с $s=1$ и
первым неравенством в (42). Тогда при соответствующем выборе
$\varphi$ и $\zeta_k=\zeta_k(1,\varphi)$ получается (41).
$\square$

\bigskip

\bigskip

Владимир Ильич Данченко, профессор,
кафедра <<Функциональный анализ и его приложения>>

Владимирский государственный университет им. А.Г. и Н.Г.
Столетовых

Email: vdanch2012@yandex.ru

\bigskip

Лев Андреевич Семин, аспирант, кафедра <<Физики и прикладной
математики>>

Владимирский государственный университет им. А.Г. и
Н.Г.Столетовых

Email: semin.lev@gmail.com

{600000, Владимир, ул. Горького, д. 87.}

\end{document}